\documentclass[11pt]{article}
\usepackage{mathrsfs}
\usepackage[justification=centering]{caption}
\usepackage{booktabs}
\usepackage[numbers,sort&compress]{natbib}
\usepackage[colorlinks,
            linkcolor=blue,
            anchorcolor=green,
            citecolor=magenta
            ]{hyperref}
\usepackage{mathrsfs,amssymb,amsfonts}
\usepackage{graphicx}
\usepackage{amsmath,amssymb,latexsym,color}
\usepackage[mathscr]{eucal}
\usepackage[numbers]{natbib}
\newtheorem{thm}{Theorem}[section]

\newtheorem{prop}[thm]{Proposition}
\newtheorem{lem}[thm]{Lemma}

\newtheorem{example}{Example}
\newtheorem{remark}{Remark}[section]
\newtheorem{pro}{Problem}[section]
\newtheorem{observation}[thm]{Observation}

\newcommand{\proof}{{\it Proof.\quad}}
\newcommand{\qed}{\hfill\Box\medskip}
\usepackage{CJK}
\begin{document}
\begin{CJK*}{GBK}{song}
\renewcommand{\abovewithdelims}[2]{
\genfrac{[}{]}{0pt}{}{#1}{#2}}

\title{\bf Strong rainbow connection numbers
of  toroidal meshes}

\author{Yulong Wei\footnote{\scriptsize Corresponding author. \newline
 {\em E-mail address:} yulong.wei@mail.bnu.edu.cn (Y. Wei), xum@bnu.edu.cn (M. Xu), wangks@bnu.edu.cn (K. Wang).}\quad Min Xu\quad Kaishun Wang\\
{\footnotesize   \em  Sch. Math. Sci. {\rm \&} Lab. Math. Com. Sys., Beijing Normal University, Beijing, 100875,  China} }
 \date{}
 \date{}
 \maketitle

\begin{abstract}

In 2011, Li et al. \cite{LLL} obtained an upper bound of the strong rainbow connection number of an $r$-dimensional undirected toroidal mesh. In this paper, this bound is improved. As a result, we give a negative answer to their  problem.

\medskip
\noindent {\em Key words:} toroidal mesh; (strong) rainbow path; (strong) rainbow connection number; Cayley graph.

\end{abstract}

\section{Introduction}

All graphs considered in this paper are finite, connected and simple. We refer to the
book \cite{BMG} for graph theory notation and terminology not described here. Let $\Gamma$ be a graph. Denote by $V(\Gamma)$ and $E(\Gamma)$ the vertex set and edge set of $\Gamma$, respectively. A sequence of distinct vertices is a path if any two consecutive vertices are adjacent. A path $P:(v_1, v_2, \ldots, v_k)$ is a cycle if $v_1$ is adjacent to $v_k$, denoted by $C_k$. The distance, $d(u, v)$, between vertices $u$ and $v$ is equal to the length of a shortest path connecting $u$ and $v$. The diameter of $\Gamma$, $d(\Gamma)$, is the maximum distance between two vertices in $\Gamma$ over all pairs of vertices.

Define a $k$-edge-coloring $\zeta: E(\Gamma)\rightarrow \{1,2,\ldots,k\}$, $k\in \mathbb{N}$,
where adjacent edges may be colored the same. A path is {\em rainbow} if no two edges of
it are colored the same. A path from $u$ to $v$ is called a {\em strong rainbow path} if it's a rainbow path with length $d(u, v)$.
If any two distinct vertices $u$ and $v$ of $\Gamma$ are connected by a (strong) rainbow path, then $\Gamma$ is called {\em (strong) rainbow-connected} under the coloring $\zeta$, and $\zeta$ is called a {\em (strong) rainbow $k$-coloring} of $\Gamma$. The {\em (strong) rainbow connection number} of $\Gamma$, denoted by (src$(\Gamma)$) rc$(\Gamma)$, is the minimum $k$ for which there exists a (strong) rainbow $k$-coloring of $\Gamma$. Clearly, we have $d(\Gamma)\leq$ rc$(\Gamma)\leq$ src$(\Gamma)$.

The (strong) rainbow connection number of a graph was first introduced by Chartrand et al. \cite{CJMZ}. Ananth and Nasre \cite{AN} proved that, for every integer $k\geq3$, deciding whether ${\rm src}(\Gamma)\leq k$  is NP-hard even when $\Gamma$ is bipartite.
(Strong) rainbow connection numbers of some special graphs have been studied in the literature, such as outerplanar graphs \cite{DXX}, Cayley graphs \cite{LLL,ML}, line graphs \cite{LSue}, power graphs \cite{MFW}, undirected double-loop networks \cite{SY} and non-commuting graphs \cite{WMW}.

The Cartesian product of two graphs $\Gamma$ and $\Lambda$ is the graph $\Gamma  \square \Lambda$ whose vertex
set is the set $\{\gamma\lambda | \gamma\in V(\Gamma), \lambda\in V(\Lambda)\}$, and two vertices $\gamma\lambda, \gamma' \lambda'$ are adjacent if $\lambda=\lambda'$ and $\{\gamma,\gamma'\}\in E(\Gamma)$ or if $\gamma=\gamma'$ and $\{\lambda,\lambda'\}\in E(\Lambda)$. The Cartesian product operation is commutative and associative, hence the Cartesian product of more factors is well-defined. The graph $C_{n_1}\square\cdots\square C_{n_r}$ is an $r$-dimensional undirected toroidal mesh, where $n_k\geq2$ for $1\leq k\leq r$.

In 2011, Li et al. proved the following theorem and proposed an open problem.

\begin{thm}{\rm\citep[Corollary 2]{LLL}}\label{HL}
Let $C_{n_k}$, $n_k\geq2$, $1\leq k \leq r$ be cycles. Then $$\sum_{1\leq k \leq r}\left\lfloor\frac{n_k}{2}\right\rfloor\leq {\rm rc}(C_{n_1}\square\cdots\square C_{n_r})\leq {\rm src}(C_{n_1}\square\cdots\square C_{n_r})\leq\sum_{1\leq k \leq r}\left\lceil\frac{n_k}{2}\right\rceil.$$ Moreover, if $n_k$ is even for every $1\leq k \leq r$, then $${\rm rc}(C_{n_1}\square\cdots\square C_{n_r})={\rm src}(C_{n_1}\square\cdots\square C_{n_r})=\sum_{1\leq k \leq r}\frac{n_k}{2}.$$
\end{thm}
\begin{pro}{\rm\citep[Remark 2]{LLL}}\label{pro1}
Given an Abelian group $G$ and an inverse closed minimal generating set $S\subseteq G \setminus 1$ of $G$, is it true that $${\rm src}(C(G,S))={\rm rc}(C(G,S))=\sum_{a\in S^*}\left\lceil\frac{|a|}{2}\right\rceil?$$
where $S^*\subseteq S$ is a minimal generating set of $G$.
\end{pro}

In this paper, we improve the upper bound of ${\rm src}(C_{n_1}\square\cdots\square C_{n_r})$ in Theorem \ref{HL}. Our main result is listed below.

\begin{thm}\label{main2}
Let $C_{n_k}$, $n_k\geq2$, $1\leq k \leq r$ be cycles. Then
\begin{equation}\label{ttt}
  {\rm src}(C_{n_1}\square\cdots\square C_{n_r})\leq
\left\{
  \begin{array}{lll}
   \left\lceil\dfrac{n_1+\cdots+n_r-\mu}{2}\right\rceil, & \hbox{$0\leq \mu\leq \left\lfloor\dfrac{r}{2}\right\rfloor$}; \\
   \\
    \left\lceil\dfrac{n_1+\cdots+n_r-r+\mu}{2}\right\rceil, & \hbox{$\left\lfloor\dfrac{r}{2}\right\rfloor+1\leq \mu\leq r$,}
  \end{array}
\right.
\end{equation}
where $\mu$ is the number of even numbers among $n_1, \ldots, n_r$.
\end{thm}

Note that an $r$-dimensional undirected toroidal mesh is a Cayley graph. As a result, Theorem \ref{main2} gives a negative answer to Problem \ref{pro1}.

\section{Preliminary results}
In this section, we will introduce some useful results for the strong rainbow connection numbers of graphs.
\begin{lem}{\rm\citep[Proposition 2.1]{CJMZ}}\label{cond}
For each integer $n\geq4$, ${\rm rc}(C_n)={\rm src}(C_n)=\left\lceil\frac{n}{2}\right\rceil$.
\end{lem}

We make the following simple observation, which we will use repeatedly.
\begin{observation}\label{dg}
Let $\Gamma$ and $\Lambda$ be two connected graphs. Then $${\rm src}(\Gamma) \leq {\rm src}(\Gamma\square \Lambda)\leq {\rm src}(\Gamma)+{\rm src}(\Lambda). $$
\end{observation}

\begin{lem}\label{zdj}
For each integer $n\geq3$, ${\rm src}(C_n\square C_2)=\left\lceil\frac{n+1}{2}\right\rceil$.
\end{lem}
\proof
Write $(0, 1, \ldots, n-1)$ for $C_n$ and $(0, 1)$ for $C_2$. Since the diameter of $C_n\square C_2$ is $\left\lceil\frac{n+1}{2}\right\rceil$, it suffices to show that ${\rm src}(C_n\square C_2)\leq\left\lceil\frac{n+1}{2}\right\rceil$. We only need to construct a strong rainbow $\left\lceil\frac{n+1}{2}\right\rceil$-coloring. Now we divide our discussion into two cases.

\medskip

\noindent {\em Case 1~~ $n=2k$}.

Define an edge-coloring $f_1$ of the graph $C_{2k}\square C_2$ by
\begin{equation*}\label{}
f_1(e)=
\left\{
  \begin{array}{lll}
   k, & \hbox{if $e=\{i0, i1\}$}; \\
    i, & \hbox{if $e=\{ij, (i+1)j\}, 0\leq i\leq k-1$};\\
     i-k, & \hbox{if $e=\{ij, (i+1)j\}$, $k\leq i\leq 2k-2$};\\
     k-1, & \hbox{if $e=\{(2k-1)j, 0j\}$}.
  \end{array}
\right.
\end{equation*}
\begin{figure}[hptb]
  \centering
  \includegraphics[width=8cm]{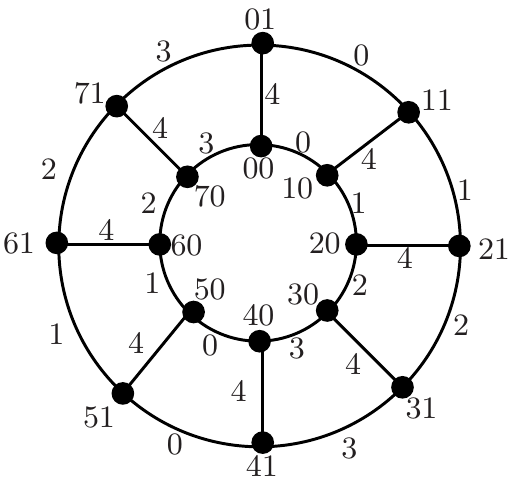}\\
  \caption{A strong rainbow $5$-coloring of $C_8\square C_2$}\label{f1}
\end{figure}
For illustration, we give a strong rainbow $5$-coloring of $C_8\square C_2$ in Figure \ref{f1}.

Note that any path from $u$ to $v$ in Table \ref{tab:1} is a strong rainbow path under the coloring $f_1$. It follows that $f_1$ is a strong rainbow $(k+1)$-coloring.

\begin{table}[htbp]\tiny
 \caption{\quad \label{tab:1} A path from $u$ to $v$ in $C_{2k}\square C_2$}
 \begin{tabular}{|c|c|c|l|}
  \toprule
  $u$ & $v$ & Condition & A path from $u$ to $v$\\
  \midrule
  $is$ & $js$ &  ~~~~~~~~~~~~~~~~~~~~~~$1\leq j-i\leq k$ ~~~~~~~~~~~~~~~~~~~~~~& $(is, (i+1)s, \ldots, js)$~~~~~~~~~~~~~~~~~~~~~~~~~~~~~~~~~~~~~~~~~\\
  $is$ & $js$ &  $ j-i\geq k+1$ & $(js, (j+1)s, \ldots, (2k-1)s, 0s, 1s, \ldots, is)$\\
  $is$ & $js$ &  $1\leq i-j\leq k$ & $(js, (j+1)s, \ldots, is)$\\
  $is$ & $js$ &  $ i-j\geq k+1$ & $(is, (i+1)s, \ldots, (2k-1)s, 0s, 1s, \ldots, js)$\\
  $i0$ & $j1$ &  $1\leq j-i\leq k$ & $(i0, i1, (i+1)1, \ldots, j1)$\\
  $i0$ & $j1$ &  $ j-i\geq k+1$ & $(j0, j1, (j+1)1, \ldots, (2k-1)1, 01, 11, \ldots, i1)$\\
  $i0$ & $j1$ &  $1\leq i-j\leq k$ & $(j0, j1, (j+1)1, \ldots, i1)$\\
  $i0$ & $j1$ &  $ i-j\geq k+1$ & $(i0, i1, (i+1)1, \ldots, (2k-1)1, 01, 11, \ldots, j1)$\\

  \bottomrule
 \end{tabular}
\end{table}

\medskip

\noindent {\em Case 2~~ $n=2k+1$}.

Define an edge-coloring $f_2$ of the graph $C_{2k+1}\square C_2$ by
\begin{equation*}\label{}
f_2(e)=
\left\{
  \begin{array}{lll}
   i, & \hbox{if $e=\{i0, i1\}, 1\leq i\leq k$}; \\
    i-k-1, & \hbox{if $e=\{i0, i1\}, k+1\leq i\leq 2k$};\\
    k, & \hbox{if $e=\{00, 01\}$};\\
    i, & \hbox{if $e=\{ij, (i+1)j\}$, $0\leq i\leq k$};\\
    i-k-1, & \hbox{if $e=\{ij, (i+1)j\}$, $k+1\leq i\leq 2k-1$};\\
    k-1, & \hbox{if $e=\{(2k)j, 0j\}$}.
  \end{array}
\right.
\end{equation*}
\begin{figure}[hptb]
  \centering
  \includegraphics[width=8cm]{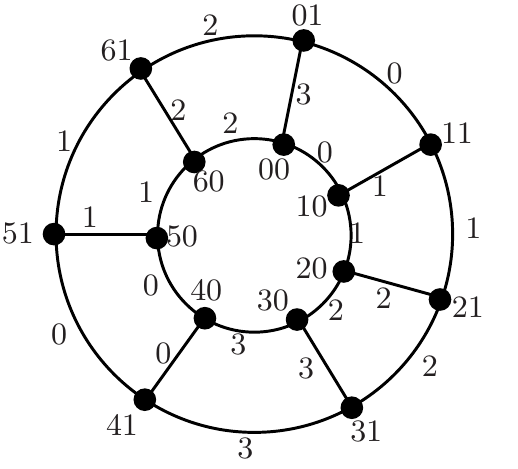}\\
  \caption{A strong rainbow $4$-coloring of $C_7\square C_2$}\label{f2}
\end{figure}
For illustration, we give a strong rainbow $4$-coloring of $C_7\square C_2$ in Figure \ref{f2}.

Note that any path from $u$ to $v$ in Table \ref{tab:2} is a strong rainbow path under the coloring $f_2$. It follows that $f_2$ is a strong rainbow $(k+1)$-coloring. $\qed$

\begin{table}[htbp]\tiny
 \caption{\quad \label{tab:2} A path from $u$ to $v$ in $C_{2k+1}\square C_2$}
 \begin{tabular}{|c|c|c|l|}
  \toprule
  $u$ & $v$ & Condition & A path from $u$ to $v$\\
  \midrule
  $is$ & $js$ &  $1\leq j-i\leq k$ & $(is, (i+1)s, \ldots, js)$\\
  $is$ & $js$ &  $ j-i\geq k+1$ & $(js, (j+1)s, \ldots, (2k)s, 0s, 1s, \ldots, is)$\\
  $is$ & $js$ &  $1\leq i-j\leq k$ & $(js, (j+1)s, \ldots, is)$\\
  $is$ & $js$ &  $ i-j\geq k+1$ & $(is, (i+1)s, \ldots, (2k)s, 0s, 1s, \ldots, js)$\\
  $00$ & $j1$ & $0\leq j\leq k$ & $(00, 01, 11, \ldots, j1)$\\
  $00$ & $j1$ & $k+1\leq j\leq 2k$ & $(00, 01, (2k)1, \ldots, j1)$\\
  $10$ & $j1$ & $1\leq j\leq k+1$ & $(10, 20, \ldots, j0, j1)$\\
  $10$ & $j1$ & $j=0$ or $k+2\leq j\leq 2k$ & $(10, 00, 01, (2k)1, \ldots, j1)$\\
  $i0$ & $j1$ & $2\leq i\leq k-1$ and $0\leq j\leq i$ & $(i0, i1, (i-1)1, \ldots, j1)$\\
  $i0$ & $j1$ & $2\leq i\leq k-1$ and $i+1\leq j\leq i+k$ & $(i0, (i+1)0, \ldots, j0, j1)$\\
  $i0$ & $j1$ & $2\leq i\leq k-1$ and $i+k+1\leq j\leq 2k$ & $(i0, (i-1)0, \ldots, 00, 01, (2k)1, (2k-1)1, \ldots, j1)$\\
  $k0$ & $j1$ & $0\leq j\leq k$ & $(k0, k1, (k-1)1, \ldots, j1)$\\
  $k0$ & $j1$ & $k+1\leq j\leq 2k$ & $(k0, (k+1)0, \ldots, j0, j1)$\\
  $i0$ & $j1$ & $k+1\leq i\leq 2k-1$ and  $0\leq j\leq i-k-1$ & $(i0, (i+1)0, \ldots, 00, 01, \ldots, j1)$\\
  $i0$ & $j1$ & $k+1\leq i\leq 2k-1$ and  $i-k\leq j\leq i$ & $(i0, i1, (i-1)1, \ldots, j1)$\\
  $i0$ & $j1$ & $k+1\leq i\leq 2k-1$ and $i+1\leq j\leq 2k$ & $(i0, (i+1)0, \ldots, j0, j1)$\\
  $(2k)0$ & $j1$ & $0\leq j\leq k-1$ & $((2k)0, 00, 01, \ldots, j1)$\\
  $(2k)0$ & $j1$ & $k\leq j\leq 2k$ & $((2k)0, (2k)1, (2k-1)1, \ldots, j1)$\\

  \bottomrule
 \end{tabular}
\end{table}

In the graph $\Gamma\square \Lambda$, we write $\Gamma y$ for $\Gamma\square\{y\}$, where $y\in V(\Lambda)$.
The {\em union} of graphs $\Gamma$ and $\Lambda$ is the graph $\Gamma\cup \Lambda$ with vertex set $V(\Gamma)\cup V(\Lambda)$ and edge set $E(\Gamma)\cup E(\Lambda)$.
\begin{prop}\label{s6}
Let $\Gamma$ be a connected graph. Then
\begin{equation}\label{ab}
  {\rm src}(\Gamma\square C_n)\leq \left\lceil\frac{n-2}{2}\right\rceil+{\rm src}(\Gamma\square C_2), n\geq3.
\end{equation}
\end{prop}
\proof
Write $(x_1, x_2, \ldots, x_n)$ for $C_n$. Meanwhile, we use $E[\Gamma x_i, \Gamma x_j]$ to denote the edge set between $\Gamma x_i$ and $\Gamma x_j$.

By Observation \ref{dg}, we have $${\rm src}(\Gamma\square C_3)\leq {\rm src}(C_3)+{\rm src}(\Gamma)=1+{\rm src}(\Gamma)\leq 1+{\rm src}(\Gamma\square C_2). $$ Thus, (\ref{ab}) holds for $n=3$. Now, we suppose that $n\geq4$.

\begin{figure}[hbt]
  \centering
  \includegraphics[width=12cm]{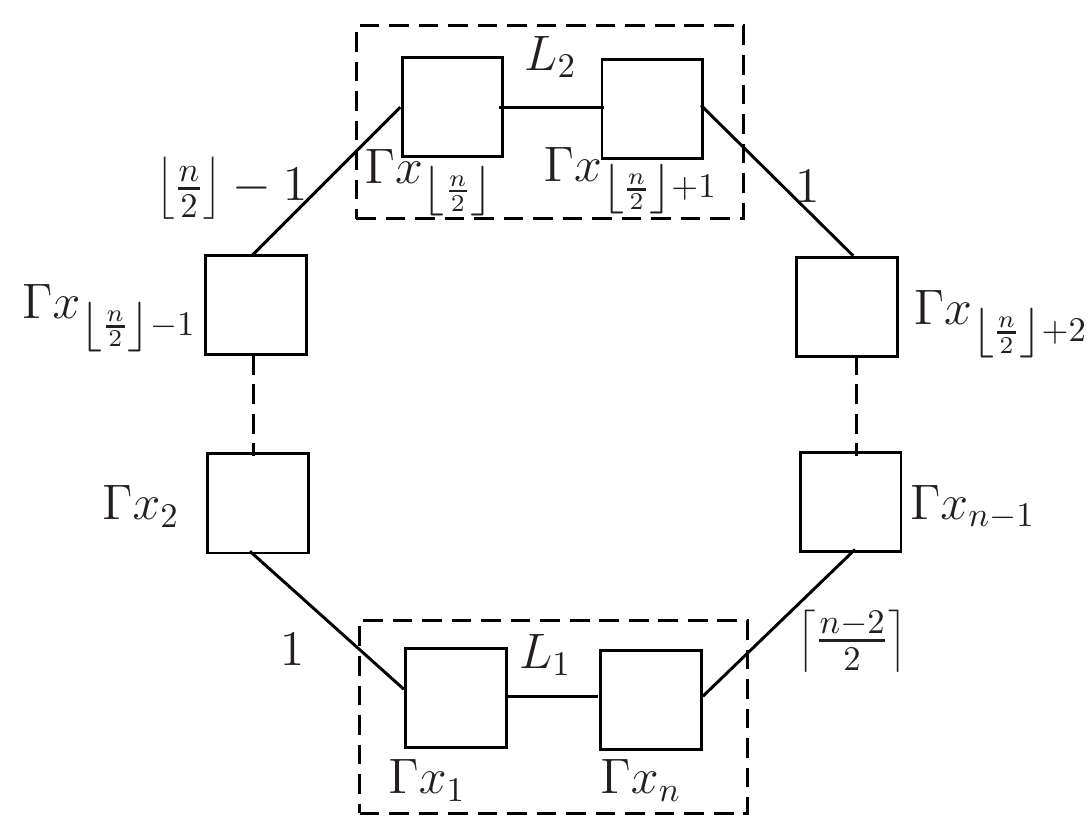}\\
  \caption{Illustration of $f_3$}\label{f6}
\end{figure}

Let $L_1$ and $L_2$ be induced subgraphs of $\Gamma\square C_n$  whose vertex sets are $V(\Gamma x_1)\cup V(\Gamma x_n)$ and $V(\Gamma x_{\left\lfloor\frac{n}{2}\right\rfloor})\cup V(\Gamma x_{\left\lfloor\frac{n}{2}\right\rfloor+1})$ respectively. Then each $L_i$ is isomorphic to $\Gamma\square C_2$. Let $S_0=\{1, \ldots, \left\lceil\frac{n-2}{2}\right\rceil\}$. Suppose that $f_1: E(\Gamma)\rightarrow S_1$ is a strong rainbow ${\rm src}(\Gamma)$-coloring of $\Gamma$, and $f_{2, i}: E(L_i)\rightarrow S_2$ is a strong rainbow ${\rm src}(\Gamma\square C_2)$-coloring of $L_i$ for $1\leq i\leq2$, where $S_0\cap S_2=\emptyset$. By Observation \ref{dg}, we may assume that $S_1\subseteq S_2$.
Define an edge-coloring $f_3: E(\Gamma\square C_n)\rightarrow S_0\cup S_2$ by
\begin{equation*}\label{}
f_3(e)=
\left\{
  \begin{array}{lll}
   i, & \hbox{if $e\in E[\Gamma x_i, \Gamma x_{i+1}]$}, \hbox{$1\leq i\leq\left\lfloor\frac{n}{2}\right\rfloor-1$}; \\
   i-\left\lfloor\frac{n}{2}\right\rfloor, & \hbox{if $e\in E[\Gamma x_i, \Gamma x_{i+1}]$}, \hbox{$\left\lfloor\frac{n}{2}\right\rfloor+1\leq i\leq n-1$}; \\
   f_{2,i}(e) , & \hbox{if $e\in E(L_i)$ for $1\leq i\leq2$};\\
    f_1(\{y_1, y_2\}), & \hbox{if $e=\{y_1x_i, y_2x_i\}\in E[\Gamma x_i, \Gamma x_i]$, $n\geq 5$ and $i\in I$},
  \end{array}
\right.
\end{equation*}
where $I=\{2, 3, \ldots, \left\lfloor\frac{n}{2}\right\rfloor-1, \left\lfloor\frac{n}{2}\right\rfloor+2, \ldots, \left\lfloor\frac{n}{2}\right\rfloor+3, \ldots, n-1\}$.  For illustration of $f_3$, see Figure \ref{f6}.

Pick any two distinct vertices $u$ and $v$ of $\Gamma\square C_n$. Write $u=y_1x_i$ and $v=y_2x_j$. Without loss of generality, we may assume that $i\leq j$.  We only need to show that there exists a strong rainbow path from $u$ to $v$ under $f_3$.
 If $i=j$, the desired result is obvious. Assume that $i\neq j$. We divide our discussion into three cases.

\noindent {\em Case 1}~~ $1\leq j-i\leq\left\lfloor\frac{n}{2}\right\rfloor$, and $2\leq j\leq\left\lfloor\frac{n}{2}\right\rfloor$ or $\left\lfloor\frac{n}{2}\right\rfloor+1\leq i\leq n-1$.

Pick a strong rainbow path $P_1$ from $y_1x_j$ to $v$ in $\Gamma x_j$. Then $$(u=y_1x_i, y_1x_{i+1}, \ldots, y_1x_j)\cup P_1$$ is a desired strong rainbow path.

\noindent {\em Case 2}~~ $1\leq j-i\leq\left\lfloor\frac{n}{2}\right\rfloor$, $j\geq\left\lfloor\frac{n}{2}\right\rfloor+1$ and $i\leq\left\lfloor\frac{n}{2}\right\rfloor$.

Pick a strong rainbow path $P_2$ from $y_1x_{\left\lfloor\frac{n}{2}\right\rfloor}$ to $y_2x_{\left\lfloor\frac{n}{2}\right\rfloor+1}$ in $L_2$. Then $$(u=y_1x_i, y_1x_{i+1}, \ldots, y_1x_{\left\lfloor\frac{n}{2}\right\rfloor}) \cup P_2 \cup(y_2x_{\left\lfloor\frac{n}{2}\right\rfloor+1}, y_2x_{\left\lfloor\frac{n}{2}\right\rfloor+2}, \ldots, y_2x_j=v)$$ is a desired strong rainbow path.

\noindent {\em Case 3}~~ $j-i\geq\left\lfloor\frac{n}{2}\right\rfloor+1$.

Pick a strong rainbow path $P_3$ from $y_1x_1$ to $y_2x_n$ in $L_1$. Then $$(u=y_1x_i, y_1x_{i-1}, \ldots, y_1x_1)\cup P_3 \cup(y_2x_n, y_2x_{n-1},\ldots, y_2x_{j+1}, y_2x_j=v)$$ is a desired strong rainbow path.

As mentioned above, we obtain the desired result. $\qed$

\section{Proof of Theorem \ref{main2}}\label{rcd}
\begin{prop}\label{s7}
Let $C_{n_k}$, $n_k\geq2$, $1\leq k \leq r$ be cycles. Then
\begin{equation}\label{cd}
{\rm src}(C_{n_1}\square\cdots\square C_{n_r})\leq\left\lceil\frac{n_1+\cdots+n_r}{2}\right\rceil.
\end{equation}
\end{prop}
\proof
Without loss of generality, we may assume that $n_1\geq n_2\geq\cdots\geq n_r$. We distinguish two cases.

\medskip

\noindent {\em Case 1~~ $n_r\geq3$}.

We prove this proposition by induction on $r$. If $r=1$, (\ref{cd}) is derived from Lemma \ref{cond}. Suppose $r=2$. If $n_2$ is even, then by Lemma \ref{cond} and Observation \ref{dg}, (\ref{cd}) holds.
If $n_2$ is odd, then by Proposition \ref{s6} and Lemma \ref{zdj}, we have
\begin{eqnarray*}
{\rm src}(C_{n_1}\square C_{n_2})&\le & \left\lceil\frac{n_2-2}{2}\right\rceil+{\rm src}(C_{n_1}\square C_2)=\left\lceil\frac{n_1+n_2}{2}\right\rceil.
\end{eqnarray*}

Now, Suppose $r\geq3$.

If each $n_i$ is odd, then
\begin{eqnarray*}
&&{\rm src}(C_{n_1}\square\cdots\square C_{n_{r-1}}\square C_{n_r})\\
&\leq &\left\lceil\frac{n_r-2}{2}\right\rceil+{\rm src}((C_{n_1}\square\cdots\square C_{n_{r-2}})\square (C_{n_{r-1}}\square C_2))~~({\rm by~Proposition}~\ref{s6})\\
&\leq&\left\lceil\frac{n_r-2}{2}\right\rceil+{\rm src}(C_{n_1}\square\cdots\square C_{n_{r-2}})+{\rm src}(C_{n_{r-1}}\square C_2)~~({\rm by~Observation}~\ref{dg})\\
&\leq&\left\lceil\frac{n_r-2}{2}\right\rceil+\left\lceil\frac{n_1+
\cdots+n_{r-2}}{2}\right\rceil+\left\lceil\frac{n_{r-1}+1}{2}\right\rceil~~({\rm by~induction~hypothesis}~{\rm and~Lemma}~\ref{zdj})\\
&=&\left\lceil\frac{n_1+\cdots+n_r}{2}\right\rceil.
\end{eqnarray*}

If $n_i$ is even for some $i$, then
\begin{eqnarray*}
&&{\rm src}(C_{n_1}\square\cdots\square C_{n_r})\\
&\leq& {\rm src}(C_{n_1}\square\cdots \square C_{n_{i-1}}\square C_{n_{i+1}}\square\cdots\square C_{n_r})+{\rm src}(C_{n_i})~~({\rm by~Observation}~\ref{dg})\\
&\leq& \left\lceil\frac{n_1+\cdots+n_{i-1}+n_{i+1}+\cdots+n_r}{2}\right\rceil+\left\lceil\frac{n_i}{2}\right\rceil~~({\rm by~induction~hypothesis})\\
&=&\left\lceil\frac{n_1+\cdots+n_r}{2}\right\rceil.
\end{eqnarray*}

\medskip

\noindent {\em Case 2~~ $n_r=2$}.

Suppose $s$ is the minimum positive integer such that $n_s=2$.
\medskip

\noindent {\em Case 2.1~~ $s=1$}.
By Observation \ref{dg}, (\ref{cd}) is obtained.

\medskip

\noindent {\em Case 2.2~~ $s\geq 2$}.
In this case, we have
\begin{eqnarray*}
&&{\rm src}(C_{n_1}\square\cdots\square C_{n_r})\\
&\leq&{\rm src}(C_{n_1}\square\cdots\square C_{n_{s-1}})+{\rm src}(C_2\square\cdots\square C_2)~~({\rm by~Observation}~\ref{dg})\\
&\leq& \left\lceil\frac{n_1+\cdots+n_{s-1}}{2}\right\rceil+r-s+1~~({\rm by~Case}~1)\\
&=&\left\lceil\frac{n_1+\cdots+n_r}{2}\right\rceil.
\end{eqnarray*}

Combining Case $1$ and Case $2$, we obtain the desired result. $\qed$

\bigskip

\noindent{\em Proof of Theorem \ref{main2}: }
If $r=1$, (\ref{ttt}) is obvious. Now, suppose $r\geq 2$. We divide our discussion into three cases.
\bigskip

\noindent {\em Case 1~~ $r=2$}.
\medskip

If $n_1=n_2=2$, Observation \ref{dg} implies that (\ref{ttt}) holds. Now suppose $n_j\geq 3$, for some $j$. Without loss of generality, we may assume that $n_1\geq 3$. By Proposition \ref{s6} and Lemma \ref{zdj}, we have
\begin{eqnarray}\label{en}
{\rm src}(C_{n_1}\square C_{n_2})\leq \left\lceil\frac{n_2-2}{2}\right\rceil+{\rm src}(C_{n_1}\square C_2)
=\left\lceil\frac{n_2-2}{2}\right\rceil+\left\lceil\frac{n_1+1}{2}\right\rceil.
\end{eqnarray}

If $\mu=1$ and $n_1$ is even, then
$${\rm src}(C_{n_1}\square C_{n_2})={\rm src}(C_{n_2}\square C_{n_1})\leq\left\lceil\frac{n_1-2}{2}\right\rceil+{\rm src}(C_{n_2}\square C_2)=\left\lceil\frac{n_1+n_2-1}{2}\right\rceil. $$ Otherwise, (\ref{en}) implies (\ref{ttt}).

\medskip

\noindent {\em Case 2~~ $r=3$}.

\medskip

If $\mu=0$ or $\mu=3$, (\ref{ttt}) holds by Proposition \ref{s7}. Now suppose $\mu=1$ or $\mu=2$. Without loss of generality, we assume that $n_1$ is even and $n_3$ is odd. By Observation \ref{dg}, Lemma \ref{cond} and Case 1, we have
\begin{eqnarray*}
{\rm src}(C_{n_1}\square C_{n_2}\square C_{n_3})\leq {\rm src}(C_{n_1})+{\rm src}(C_{n_2}\square C_{n_3})
\leq\left\lceil\frac{n_1+n_2+n_3-1}{2}\right\rceil.
\end{eqnarray*}

\noindent {\em Case 3~~ $r\geq4$}.
\medskip

If $\mu=0$ or $\mu=r$, Proposition \ref{s7} implies (\ref{ttt}).
In the following, we assume that $n_1, \ldots, n_\mu$ are even and $n_{\mu+1}, \ldots, n_r$ are odd.

\noindent {\em Case 3.1~~ $1\leq \mu\leq \left\lfloor\frac{r}{2}\right\rfloor-1$}.
\begin{eqnarray*}
&&{\rm src}(C_{n_1}\square\cdots\square C_{n_r})\\
&\leq&{\rm src}(C_{n_{2\mu+1}}\square\cdots\square C_{n_r})+\sum_{1\leq s\leq \mu, s+t=2\mu+1} {\rm src}(C_{n_s}\square C_{n_t})~~({\rm by~Observation}~\ref{dg})\\
&\leq& \left\lceil\frac{n_{2\mu+1}+\cdots+n_r}{2}\right\rceil+\sum_{1\leq s\leq \mu, s+t=2\mu+1}\frac{n_s+n_t-1}{2}~~({\rm by~Proposition}~\ref{s7}~{\rm and~Case~1})\\
&=&\left\lceil\frac{n_1+\cdots+n_r-\mu}{2}\right\rceil.
\end{eqnarray*}

\noindent {\em Case 3.2~~ $\mu=\left\lfloor\frac{r}{2}\right\rfloor$}.
\begin{eqnarray*}
&&{\rm src}(C_{n_1}\square\cdots\square C_{n_r})\\
&\leq& \frac{1}{2}(1+(-1)^{r+1}){\rm src}(C_{n_r})+\sum_{1\leq s\leq \mu, s+t=2\mu+1} {\rm src}(C_{n_s}\square C_{n_t})~~({\rm by~Observation}~\ref{dg})\\
&\leq& \frac{1}{2}(1+(-1)^{r+1})\left\lceil\frac{n_r}{2}\right\rceil+\sum_{1\leq s\leq \mu, s+t=2\mu+1}\frac{n_s+n_t-1}{2}~~({\rm by~Proposition}~\ref{s7}~{\rm and~Case~1})\\
&=&\left\lceil\frac{n_1+\cdots+n_r-\mu}{2}\right\rceil.
\end{eqnarray*}

\noindent {\em Case 3.3~~ $\left\lfloor\frac{r}{2}\right\rfloor+1\leq \mu\leq r-1$}.
\begin{eqnarray*}
&&{\rm src}(C_{n_1}\square\cdots\square C_{n_r})\\
&\leq&{\rm src}(C_{n_1}\square\cdots\square C_{2\mu-r})+\sum_{\mu+1\leq t\leq r, s+t=2\mu+1} {\rm src}(C_{n_s}\square C_{n_t})~~({\rm by~Observation}~\ref{dg})\\
&\leq& \left\lceil\frac{n_1+\cdots+n_{2\mu-r}}{2}\right\rceil+\sum_{\mu+1\leq t\leq r, s+t=2\mu+1}\frac{n_s+n_t-1}{2}~~({\rm by~Proposition}~\ref{s7}~{\rm and~Case~1})\\
&=&\left\lceil\frac{n_1+\cdots+n_r-r+\mu}{2}\right\rceil.
\end{eqnarray*}

As mentioned above, we obtain the desired result. $\qed$

\section*{Acknowledgement}
M. Xu's research is supported by the National Natural Science Foundation of China (11571044, 61373021).
K. Wang's research is supported by the National Natural Science Foundation of China (11671043, 11371204).

\end{CJK*}

\end{document}